\documentclass[10pt]{article}
\usepackage{amsfonts, amssymb, amsmath, amscd}
\usepackage[dvips]{epsfig,graphics}
\usepackage{latexsym, longtable, tabularx, array}
\usepackage{epsfig, psfrag}

\newtheorem{thm}{Theorem}
\newtheorem{prop}{Proposition}
\newtheorem{lemma}{Lemma}
\newtheorem{cor}[thm]{Corollary}

\newtheorem{defn}{Definition}

\newcommand{\Z}{\mathbb{Z}}
\newcommand{\G}{G}
\newcommand{\bbT}{T}
\newcommand{\bbQ}{\mathbb{Q}}
\newcommand{\s}{\sigma}
\newcommand{\T}{\bar{\tau}}
\renewcommand{\t}{\bar{t}}

\renewcommand{\c}{\cite}
\newcommand{\pf}{{\em Proof: \quad }}
\newcommand{\done}{\hfill $\blacksquare$}
\newcommand{\kb}[1]{\ensuremath{\langle #1 \rangle}}
\newcommand{\del}{\partial}
\newcommand{\D}{\mathbb{G}}
\renewcommand{\H}{\mathbb{H}}
\newcommand{\C}{\mathcal{C}}

\begin{document}

\title{Graphs on surfaces and Khovanov homology}
\author{
Abhijit  Champanerkar
\footnote{Supported by National Science Foundation DMS-0455978.}\\
{\em{\small Department of Mathematics and Statistics, University of
  South Alabama}}\\ \\
 Ilya Kofman
\footnote{Supported by National Science Foundation DMS-0456227 and PSC-CUNY 60046-3637.}\\
{\em {\small Department of Mathematics, College of Staten Island, City University of New York}}\\ \\
 Neal Stoltzfus
\footnote{Supported by National Science Foundation DMS-0456275.}\\
{\em {\small Department of Mathematics, Louisiana State University}}}
\date{}

\maketitle

\centerline{\em\large{ In memory of Xiao-Song Lin}}

\begin{abstract}
\noindent
Oriented ribbon graphs (dessins d'enfant) are graphs embedded in
oriented surfaces.  A quasi-tree of a ribbon graph is a spanning
subgraph with one face, which is described by an ordered chord
diagram.  We show that for any link diagram $L$, there is an
associated ribbon graph whose quasi-trees correspond bijectively to
spanning trees of the graph obtained by checkerboard coloring $L$.
This correspondence preserves the bigrading used for the spanning tree
model of Khovanov homology, whose Euler characteristic is the Jones
polynomial of $L$.  Thus, Khovanov homology can be expressed in terms
of ribbon graphs, with generators given by ordered chord diagrams.
\end{abstract}

\section{Introduction}
A {\em ribbon graph} is a multi-graph (loops and multiple edges
allowed) that is embedded in a surface.  In the case when the surface
is oriented, the embedding determines a cyclic order on the edges at
every vertex, which is called an orientation for the ribbon graph.
Other terms for oriented ribbon graphs include: combinatorial maps,
fat graphs, cyclic graphs, graphs with rotation systems, and dessins
d'enfant (see \c{BR1}). In this paper, all ribbon graphs will be
oriented.

The Jones polynomial of any link can be obtained as a specialization
of the Bollob\'as-Riordan-Tutte polynomial of a ribbon graph obtained
from the link diagram \c{DFKLS1}.  The Jones polynomial also has an
expansion in terms of spanning trees of the Tait graph, obtained by
checkerboard coloring the link diagram.  Moreover, with an appropriate
bigrading, these spanning trees generate Khovanov homology, whose
bigraded Euler characteristic is the Jones polynomial \c{KH}.

We show that there is a one-to-one correspondence between spanning
trees of the Tait graph and {\em quasi-trees}, which are spanning
ribbon subgraphs with one face.  We translate the data used to define the
bigrading for spanning trees to the language of ribbon graphs using ordered
chord diagrams.  The Khovanov homology results in \c{KH} are then
expressed in terms of ribbon graphs and ordered chord diagrams.  This leads
to the question: Do any of the algebraic structures known for chord
diagrams carry over to Khovanov homology?

This project was inspired by \c{DFKLS2} and \c{Manturov}, and the
correspondence we establish implies some of their results.

\section{Quasi-trees and spanning trees}

Let $D$ be a connected link diagram.  A checkerboard coloring of $D$
determines the Tait graph $\G$.  An edge of $\G$ is positive if the
shaded regions of its endpoints are joined by $A$-smoothing the
corresponding crossing of $D$.  Otherwise, the edge is negative.  We
take either $\G$ or its planar dual so that $E_+(\G)\geq E_-(\G)$.

Let $\D$ be the all-$A$ ribbon graph of $D$ as defined in \c{DFKLS1}.
Let $V(\D)$ be the number of vertices of $\D$, which is the number of
components in the all-$A$ state of $D$.  A ribbon subgraph $\H\subset\D$ is
called a spanning subgraph if $V(\H)=V(\D)$.  Let $F(\H)$ be
the number of faces of $\H$, which is the number of complementary
regions in the orientable surface of minimal genus on which $\H$
embeds.  A quasi-tree $\bbQ$ is a spanning subgraph of $\D$ with
$F(\bbQ)=1$ (see Definition 3.1 of \c{DFKLS2}).

\begin{thm}\label{mainthm}
Quasi-trees of $\D$ are in one-one correspondence with spanning trees of $\G$:
$$ \bbQ_j \leftrightarrow \bbT_v \qquad {\rm where }\qquad v+j = \frac{V(\G) + E_+(\G) - V(\D)}{2} $$
$\bbQ_j$ denotes a quasi-tree of genus $j$, and $T_v$ denotes a spanning tree with $v$ positive edges.
\end{thm}

The proof will use the following lemma.
A state $s$ of $D$ is given by $s:{\rm Edges}(\D)\rightarrow \{A,B\}$.
Let $|s|$ denote the number of components in the corresponding smoothing of $D$.
In Section 4 of \c{DFKLS1}, the ribbon graph $\D(s)$ was defined such that $V(\D(s))=|s|$.
We now define a different correspondence between states of $D$ and ribbon graphs:

\begin{lemma}\label{Hs}
Spanning subgraphs $\H\subset\D$ are in one-one correspondence with states $s$ of $D$, such that
$s(\H)(e) = B$ iff $e\in\H$.  Thus, $F(\H)=|s|$ and $V(\H)=V(\D)$.
\end{lemma}
\pf
For any state $s$ of $D$, let $D_s$ denote the following link diagram:
$$ D_s  =
\begin{cases}
A{\rm -smoothing\ at\ } e  & {\rm if\ } s(e)=A \\
D{\rm -crossing\ at\ } e   & {\rm if\ } s(e)=B
\end{cases} $$
Let $\D_A(D_s)$ and $\D_B(D_s)$ denote the all-$A$ and all-$B$ ribbon graphs of $D_s$, respectively.
For any state $s$ of $D$, define $\H(s) = \D_A(D_s)$.
It follows that $s=s(\H(s))$ and $\H=\H(s(\H))$.

The ribbon graphs $\H(s)=\D_A(D_s)$ and $\D_B(D_s)$ are dual in the
sense of Lemma 4.1 \c{DFKLS1}.  By this duality, $F(\H(s))=
V(\D_B(D_s))=|s|$.  Also, $V(\H(s))=V(\D_A(D_s))=V(\D)$.  \done

The Jordan trail of a connected link diagram is a simple closed curve
obtained by smoothing each crossing \c[p. 2]{KauffmanFKT}.  There is a
one-one correspondence between Jordan trails of $D$ and spanning trees
of the Tait graph $G$ \c[p. 56]{KauffmanFKT}. In particular, the
Jordan trail of a spanning tree $T$ bounds a planar neighborhood of
$T$.

{\em Proof of Theorem \ref{mainthm}: \quad }
In the table below, let $\tau,t,\T,\t$ denote
a positive edge in $\bbT$, a positive edge in $\G-\bbT$, a negative
edge in $\bbT$, and a negative edge in $\G-\bbT$, respectively.  The
Jordan trail of $\bbT$ is then given by
the smoothings of $D$ shown in the second row.
By Lemma \ref{Hs}, each Jordan trail corresponds to a spanning
subgraph of $\D$ with one face, which is a quasi-tree.  Let $\bbQ$
be the quasi-tree that corresponds to $T$.  If $Q\in\bbQ$ and
$q\in(\D-\bbQ)$, then $s(\bbQ)$, given in Lemma \ref{Hs}, determines
the correspondence:
\begin{center}
\begin{tabular}{c|c|c|c}
$\tau $ & $ t$ & $\T$ & $\t$ \\
$A $ & $ B $ & $ B $ & $ A $\\
$q $ & $ Q$ & $Q$ & $q$
\end{tabular}
\end{center}

To prove the numerical claim, for any $\bbT$ in $\G$,
$v(\bbT)=\#\tau$ and $E(\bbQ)=\#t + \#\T$.
\begin{eqnarray*}
j&=& \frac{1-V(\bbQ)+E(\bbQ)}{2}=\frac{1-V(\D)+\#t + \#\T}{2}\\
v+j &=& \frac{2(\#\tau)+1-V(\D)+\#t + \#\T}{2}
= \frac{V(G)+E_+(G)-V(\D)}{2}
\end{eqnarray*}
since $\#\tau+\#\T=E(T)=V(G)-1$ and $\#\tau+\#\t=E_+(G)$.
\done

\section{Quasi-tree complex for Khovanov homology}

To construct the spanning tree chain complex in \c{KH}, every spanning
tree $\bbT$ of the Tait graph $\G$ was given a bigrading $(u(\bbT),
v(\bbT))$.  By Theorem \ref{mainthm}, the $v$-grading, which is the
number of positive edges in $\bbT$, is determined by the genus of the
corresponding quasi-tree $\bbQ$.  The $u$-grading, which was defined
using activities in the sense of Tutte, also has a quasi-tree analogue
in terms of the ordered chord diagram for $\bbQ$.

If $D$ has $n$ ordered crossings, let $\D$ be given by permutations
$(\s_0,\s_1,\s_2)$ of the set $\{1,\ldots,2n\}$, such that the $i$-th
crossing corresponds to half-edges $\{2i-1, 2i\}$, which are marked on
the components of the all-$A$ state of $D$.  To be precise, suppose at
the crossings of $D$, the strands are parallel to $y=x$ or $y=-x$,
then we require that the marks $2i-1$ and $2i$ be in the half-planes
$y>x$ and $y<x$, respectively.
For an example, see Figure \ref{trefoil} in Section \ref{example}.
We give the components of the all-$A$
state of $D$ the admissible orientation for which outer ones are
oriented counterclockwise (see \c{DFKLS1}).  In this way, every
component has a well-defined positive direction.

The orbits of $\s_0$ form the vertex set. In particular, $\s_0$
is given by noting the half-edge marks when going in the positive
direction around the components of the all-$A$ state of $D$.
The other permutations are given by
$ \s_1 = \prod_{i=1}^{n} (2i-1, 2i)$ and $\s_2 = \s_1 \circ \s_0^{-1}$

\begin{figure}
\begin{center}
\psfrag{1}{\tiny{$1$}}
\psfrag{2}{\tiny{$2$}}
\psfrag{3}{\tiny{$3$}}
\psfrag{4}{\tiny{$4$}}
\psfrag{5}{\tiny{$5$}}
\psfrag{6}{\tiny{$6$}}
\psfrag{7}{\tiny{$7$}}
\psfrag{8}{\tiny{$8$}}
\psfrag{3c}{\tiny{$3$}}
\psfrag{4c}{\tiny{$4$}}
\psfrag{7c}{\tiny{$7$}}
\psfrag{8c}{\tiny{$8$}}
\includegraphics[height=1.25in]{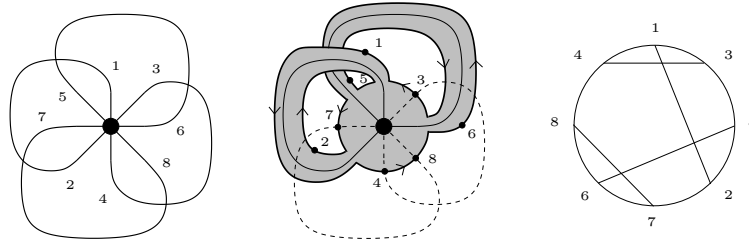}
 \end{center}
 \caption{Ribbon Graph $\D$, quasi-tree $\bbQ=(12)(56)$ with curve $\gamma_{\bbQ}$, chord diagram $C_{\bbQ}$}
 \label{cq}
 \end{figure}

Let an {\em ordered chord diagram} denote a circle marked with
$\{1,\ldots,2n\}$ in some order, and chords joining all pairs $\{2i-1, 2i\}$.
\begin{prop}\label{activity}
Every quasi-tree $\bbQ$ corresponds to the ordered chord diagram
$C_{\bbQ}$ with consecutive markings in the positive direction given
by the permutation:
$$\s(i) =
\begin{cases}
\s_0(i) & i\notin\bbQ \\
\s_2^{-1}(i) & i\in\bbQ
\end{cases} $$
\end{prop}
\pf
Since $\bbQ$ is a quasi-tree, $\gamma_{\bbQ}$ is one simple closed curve.
If we choose an orientation on $S(\D)$, we can traverse
$\gamma_{\bbQ}$ along successive boundaries of bands and vertex discs,
such that we always travel around the boundary of each disc in a
positive direction (i.e., the disc is on the left).  If a half-edge is not in $\bbQ$, $\gamma_{\bbQ}$
will pass across it travelling along the boundary of a vertex disc to
the next band. If a half-edge is in $\bbQ$, $\gamma_{\bbQ}$ traverses along
one of the edges of its band.
On $\gamma_{\bbQ}$, we mark a half-edge not in $\bbQ$ when $\gamma_{\bbQ}$
passes across it along the boundary of the vertex disc and we mark a
half-edge in $\bbQ$ when we traverse an edge of a band in the direction
of the half-edge.  If the half-edge $i$ is not in $\bbQ$, travelling
along the boundary of a vertex disc, the next half-edge is given by
$\sigma_0$.  If the half-edge $i$ is in $\bbQ$, traversing the edge of
its band to the vertex disc and then along the boundary of that disc,
the next half-edge is given by $\sigma_0\sigma_1=\sigma_2^{-1}$.  For
example, see Figure \ref{cq}.

As $\bbQ$ is a quasi-tree, each of its half-edges must be in the orbit
of its single face, while the complementary set of half-edges are met
along the boundaries of the vertex discs.
Since we mark all half-edges traversing $\gamma_{\bbQ}$, the chord diagram $C_{\bbQ}$ parametrizes $\gamma_{\bbQ}$.
\done

Note that if $\D$ is the all-A ribbon graph of a connected link diagram $D$,
by the proof of Theorem \ref{mainthm}, following $\gamma_{\bbQ}$ along an edge of
$\bbQ$ is given by the $B$-smoothing of that crossing of
$D_{s(\bbQ)}$.  Therefore, the chord diagram $C_{\bbQ}$ parametrizes
both $\gamma_{\bbQ}$ and the Jordan trail for $T$, which is the
all-$B$ state of $D_{s(\bbQ)}$.

To compute the genus $g(\bbQ)$ from $C_{\bbQ}$, let $C$ be the sub-chord diagram of
chords that correspond to edges in $\bbQ$.  Then $g(\bbQ)$ is half the rank
of the adjacency matrix of the intersection graph of $C$ \c{BR1}.

\begin{defn}
  Using $\min(i,\s_1(i))$, there is an induced total order on the
  chords of $C_{\bbQ}$. A chord is {\em live} if it does not intersect
  lower-ordered chords, and otherwise it is {\em dead}.  For any
  quasi-tree $\bbQ$, an edge $e$ is {\em live} or {\em dead} when the
  corresponding chord of $C_{\bbQ}$ is live or dead.
\end{defn}
In Figure \ref{cq}, we show $C_{\bbQ}$ such that the only edge live with respect to $\bbQ$ is $(12)$.
\begin{lemma}\label{live}
If $\bbT$ corresponds to $\bbQ$, as in Theorem \ref{mainthm}, then the
$i$-th edge of $\D$ is live with respect to $\bbQ$ if and only if the
$i$-th edge of $\G$ is live with respect to $\bbT$.
\end{lemma}
\pf
In $C_{\bbQ}$, the $i$-th and $j$-th chords intersect if and only
if going around the Jordan trail for $\bbT$ in some direction, we see
cyclic permutations of the marks $(2i-1,2j-1,2i,2j)$ or
$(2i-1,2j,2i,2j-1)$.  Now, $e_i\in\mathit{cut(\bbT,e_j)}$ or
$e_i\in\mathit{cyc(\bbT,e_j)}$ if and only if the Jordan trail becomes
disconnected when the $j$-th smoothing is changed, and is re-connected
when the $i$-th smoothing is changed.  Equivalently, $C_{\bbQ}$
becomes disconnected when unzipped along the $j$-th chord, and becomes
re-connected when unzipped along the $i$-th chord, which occurs if and
only if the $i$-th and $j$-th chords intersect:
$$\mbox{\epsfig{file=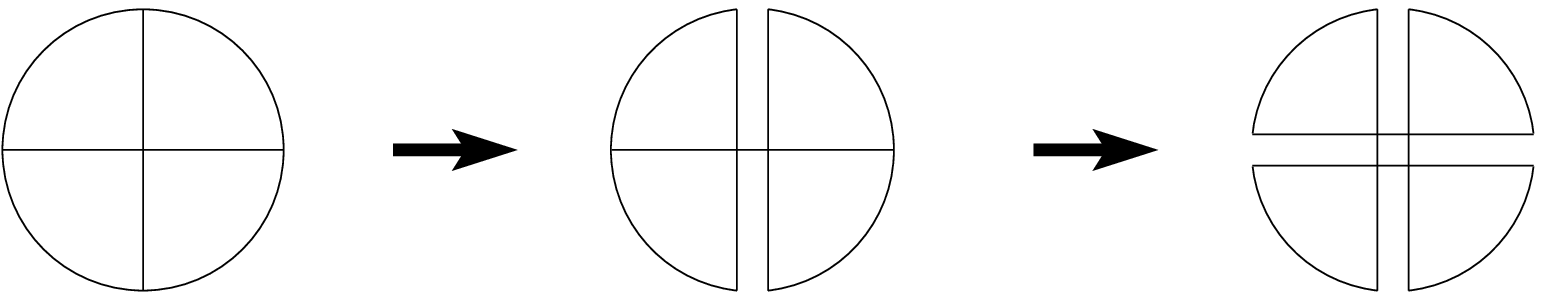,height=16pt}}$$
Therefore, $e_i$ is live with respect to $\bbT$ exactly when the
$i$-th chord does not intersect lower-ordered chords.  \done

\begin{defn}
For any quasi-tree $\bbQ$ of $\D$, we define
$$ u(\bbQ) = \#\{{\rm live \ not \ in} \ \bbQ \} - \#\{{\rm live\ in }\ \bbQ\}
\quad {\rm and }\quad v(\bbQ) = -g(\bbQ) $$
Define $\C(\D)=\oplus_{u,v} \C_v^u(\D)$, \ where \
$\C_{v}^{u}(\D) = \Z\kb{\bbQ\subset \D |\; u(\bbQ)=u,\; v(\bbQ)=v }$
\end{defn}

\begin{thm}\label{KHthm}
For a knot diagram $D$, there exists a quasi-tree complex
$\C(\D)=\{\C_v^u(\D), \del\}$ with $\del: \C_v^u \to \C_{v-1}^{u-1}$
that is a deformation retract of the reduced Khovanov complex.
In particular, the reduced Khovanov homology $\widetilde{H}^{i,j}(D;\Z)$ is given by
$$ \widetilde{H}^{i,j}(D;\Z) \cong H_v^u(\C(\D);\Z) $$
with the indices related as follows:
$$ u = j-i-w(D)+1 \quad {\rm and}\quad v = j/2 -i+(V(\D)-c_+(D))/2 $$
where $w(D)$ is the writhe, $c_+(D)$ is the number of positive
crossings of $D$, and $V(\D)$ is the number of components in the
all-$A$ state of $D$, which is the number of vertices of $\D$.
\end{thm}

\pf The result follows from Theorem 5 of \c{KH} and Theorem
\ref{mainthm}, once we establish for a spanning tree $\bbT$ corresponding
to a quasi-tree $\bbQ$ that their bigradings are related as claimed.

By Lemma \ref{live}, edges of $\bbQ$ and $\bbT$ are live exactly when
they correspond.  From \c{KH}, $u(\bbT)= \# L - \# \ell - \# \bar{L} +
\# \bar{\ell}$.  By the proof of Theorem \ref{mainthm},
$\{L,\bar{\ell}\}$-edges of $\bbT$ correspond to live edges not in
$\bbQ$, and $\{\ell,\bar{L}\}$-edges of $\bbT$ correspond to live
edges in $\bbQ$.  Therefore,
$$ u(\bbQ) =  \#\{{\rm live \ not \ in} \ \bbQ \} - \#\{{\rm live\ in } \ \bbQ\} = u(\bbT) = j-i-w(D)+1 $$

By Theorem 5 of \c{KH} and Theorem \ref{mainthm},
\begin{eqnarray*}
 v(\bbQ) & = & -g(\bbQ)  =  v(\bbT) - \frac{V(\G) + E_+(\G) - V(\D)}{2} \\
& = & \left( \frac{j}{2} -i-\frac{w(D)-k(D)-2}{4} \right) - \frac{V(\G) + E_+(\G) - V(\D)}{2} \\
& = &  \frac{j}{2} -i -\frac{w(D)+E(G)}{4} + \frac{V(\D)}{2} = \frac{j}{2} -i + \frac{V(\D)-c_+(D)}{2}
\end{eqnarray*}
where we used that $w(D)+E(\G) = w(D) + c(D) = 2c_+(D)$.
\done

The following bound for the thickness of Khovanov homology in terms of ribbon graph genus was obtained by Manturov \c{Manturov}:
\begin{cor}\label{thickness}
Let $g(\D)$ denote the genus of the all-$A$ ribbon graph of $D$.
The thickness of the reduced Khovanov homology of $D$ is less than or equal to $g(\D)+1$.
\end{cor}

\pf
For any quasi-tree $\bbQ$ of $\D$, $-g(\D)\leq v(\bbQ)\leq 0$.
Therefore, $\C(\D)$ has $g(\D)+1$ rows, so $H_v^u(\C(\D);\Z)$ has at most $g(\D)+1$ rows.
\done

Corollary \ref{thickness} is stronger than Theorem 13($ii$) of \c{KH}
because, for instance, links whose ribbon graphs have genus one are a
much richer class than $1$-almost alternating links.  (For example,
see Lemma 4.3 \c{DFKLS2}.)  However, by the correspondence in Theorem
\ref{mainthm}, the two proofs are the same: Since $0\leq j\leq g(\D)$,
$$ g(\D) = \max_{\bbT\subset \G} v(\bbT) - \min_{\bbT\subset \G} v(\bbT) $$

\section{Example}
\label{example}

\begin{figure}
\begin{center}
\psfrag{1}{\tiny{$1$}}
\psfrag{2}{\tiny{$2$}}
\psfrag{3}{\tiny{$3$}}
\psfrag{4}{\tiny{$4$}}
\psfrag{5}{\tiny{$5$}}
\psfrag{6}{\tiny{$6$}}
\psfrag{7}{\tiny{$7$}}
\psfrag{8}{\tiny{$8$}}
\psfrag{3b}{\tiny{$\overline{3}$}}
\psfrag{4b}{\tiny{$\overline{4}$}}
 \includegraphics[height=2.25in]{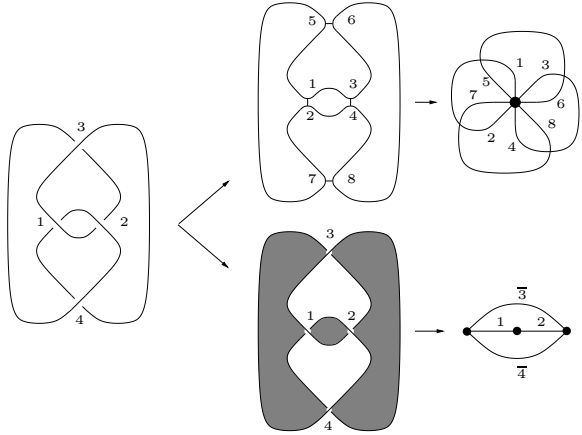}
 \end{center}
 \caption{The $4$-crossing trefoil, the all-$A$ ribbon graph $\D$, and the Tait graph $G$}
 \label{trefoil}
 \end{figure}

As an example, we use a $4$-crossing diagram of the trefoil.
In Figure \ref{trefoil}, we show the diagram $D$, the Tait graph $G$, and the all-$A$ ribbon graph $\D$, given by
\begin{center}
\begin{tabular}{m{2in} m{2in}}
\begin{tabular}{ccl}
$\sigma_0$&$=$&$(15724863)$\\
$\sigma_1$&$=$&$(12)(34)(56)(78)$\\
$\sigma_2$&$=$&$(14)(2835)(67)$
\end{tabular} &
\psfrag{1}{\tiny{$1$}}
\psfrag{2}{\tiny{$2$}}
\psfrag{3}{\tiny{$3$}}
\psfrag{4}{\tiny{$4$}}
\psfrag{5}{\tiny{$5$}}
\psfrag{6}{\tiny{$6$}}
\psfrag{7}{\tiny{$7$}}
\psfrag{8}{\tiny{$8$}}
\includegraphics[height=0.8in]{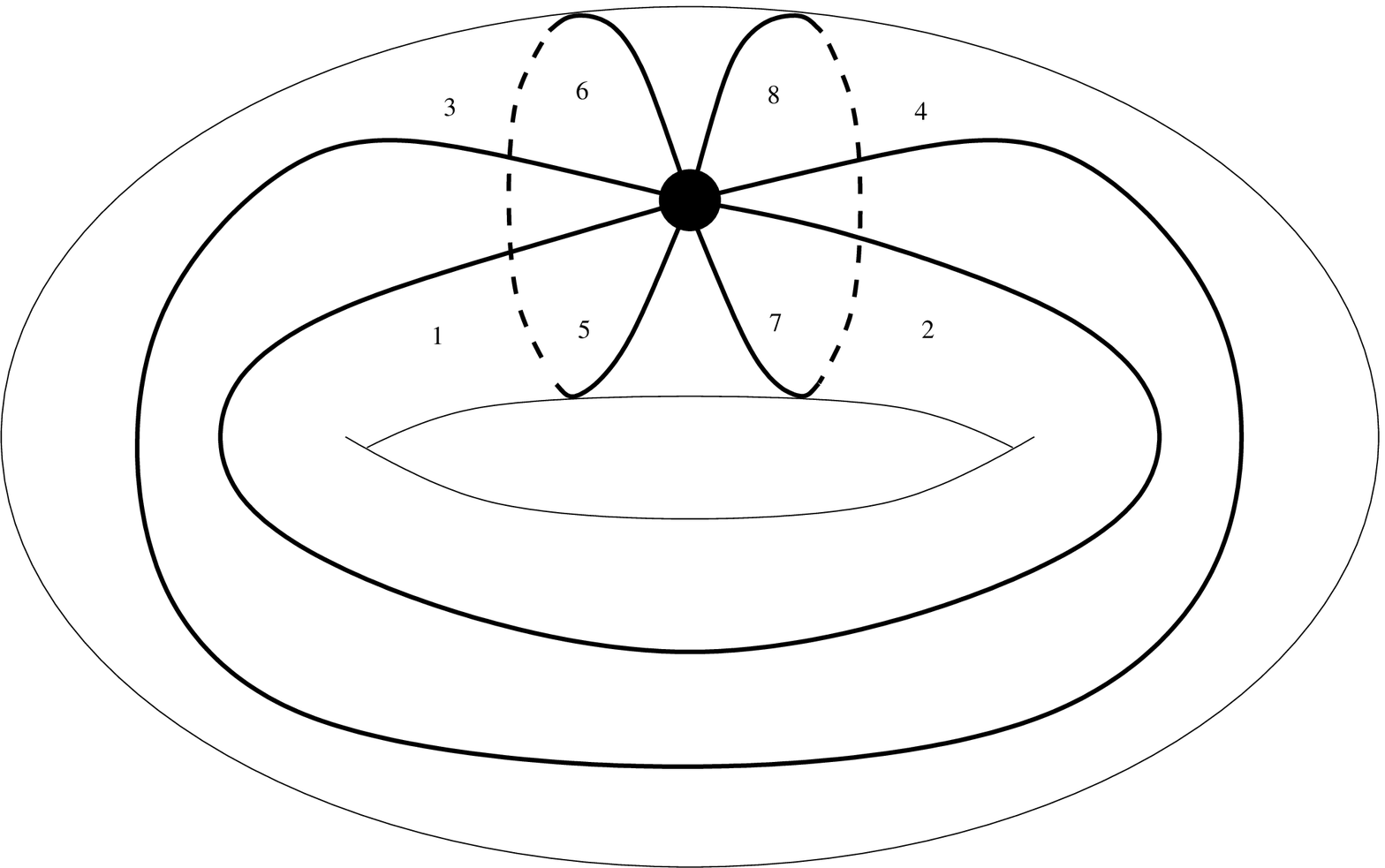}
\end{tabular}
\end{center}
\begin{tabular}{{m{3in} m{2in}}}
The ordered chord diagram for each quasi-tree is given by Corollary
\ref{activity}.  This order can be seen from the corresponding Jordan
trail, which is shown
for the quasi-tree $\bbQ_1$ with edges $(12)$ and $(56)$:
&
\makebox[2.3in]{
\psfrag{1}{\tiny{$1$}}
\psfrag{2}{\tiny{$2$}}
\psfrag{3}{\tiny{$3$}}
\psfrag{4}{\tiny{$4$}}
\psfrag{5}{\tiny{$5$}}
\psfrag{6}{\tiny{$6$}}
\psfrag{7}{\tiny{$7$}}
\psfrag{8}{\tiny{$8$}}
\psfrag{3b}{\tiny{$\overline{3}$}}
\psfrag{4b}{\tiny{$\overline{4}$}}
 \includegraphics[height=0.9in]{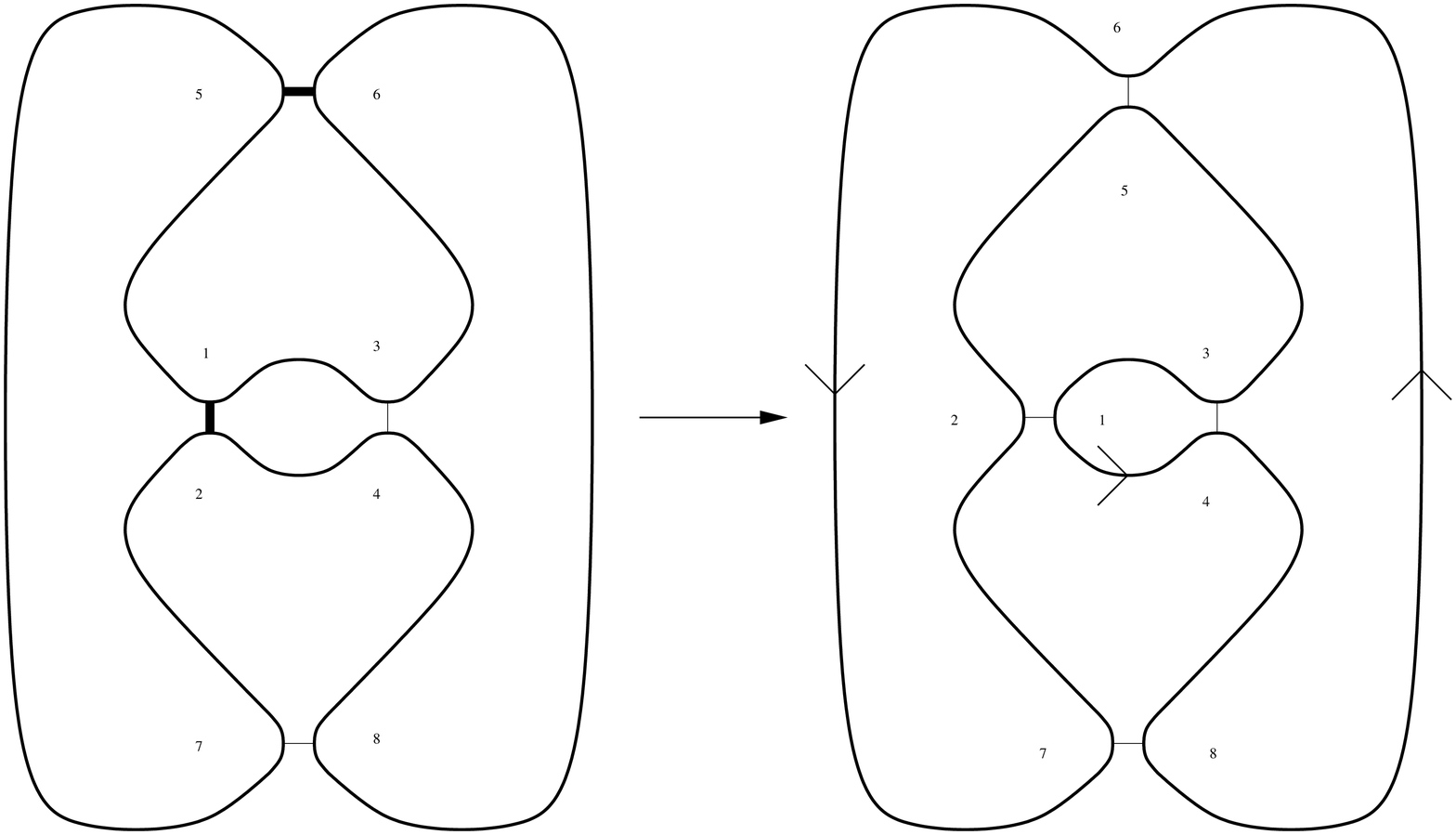}
}
\end{tabular}

Below we show the correspondence between the spanning trees of $G$,
the quasi-trees of $\D$, and their chord diagrams.  The circled
numbers on each chord diagram indicate edges in $\D-\bbQ$.  The
activities follow the convention: capital letters for edges in the
spanning tree or quasi-tree, bar for negative edges, $L$ or $\ell$ for
live, $D$ or $d$ for dead.

\begin{center}
\psfrag{1}{\tiny{$1$}}
\psfrag{2}{\tiny{$2$}}
\psfrag{3}{\tiny{$3$}}
\psfrag{4}{\tiny{$4$}}
\psfrag{5}{\tiny{$5$}}
\psfrag{6}{\tiny{$6$}}
\psfrag{7}{\tiny{$7$}}
\psfrag{8}{\tiny{$8$}}
\psfrag{3b}{\tiny{$\overline{3}$}}
\psfrag{4b}{\tiny{$\overline{4}$}}
\psfrag{1c}{\tiny{$\textcircled{1}$}}
\psfrag{2c}{\tiny{$\textcircled{2}$}}
\psfrag{3c}{\tiny{$\textcircled{3}$}}
\psfrag{4c}{\tiny{$\textcircled{4}$}}
\psfrag{5c}{\tiny{$\textcircled{5}$}}
\psfrag{6c}{\tiny{$\textcircled{6}$}}
\psfrag{7c}{\tiny{$\textcircled{7}$}}
\psfrag{8c}{\tiny{$\textcircled{8}$}}
\psfrag{h1}{\footnotesize{$\bbQ_1=LdDd$}}
\psfrag{h2}{\footnotesize{$\bbQ_2=Ld\ell D$}}
\psfrag{h3}{\footnotesize{$\bbQ_3=\ell DDd$}}
\psfrag{h4}{\footnotesize{$\bbQ_4=\ell D\ell D$}}
\psfrag{h5}{\footnotesize{$\bbQ_5=\ell \ell dd$}}
\psfrag{t1}{\footnotesize{$T_1=\ell D \overline{D}\overline{d}$}}
\psfrag{t2}{\footnotesize{$T_2=\ell D \overline{\ell}\overline{D}$}}
\psfrag{t3}{\footnotesize{$T_3=L d \overline{D}\overline{d}$}}
\psfrag{t4}{\footnotesize{$T_4=L d \overline{\ell}\overline{D}$}}
\psfrag{t5}{\footnotesize{$T_5=L L \overline{d}\overline{d}$}}
\includegraphics[width=4.25in]{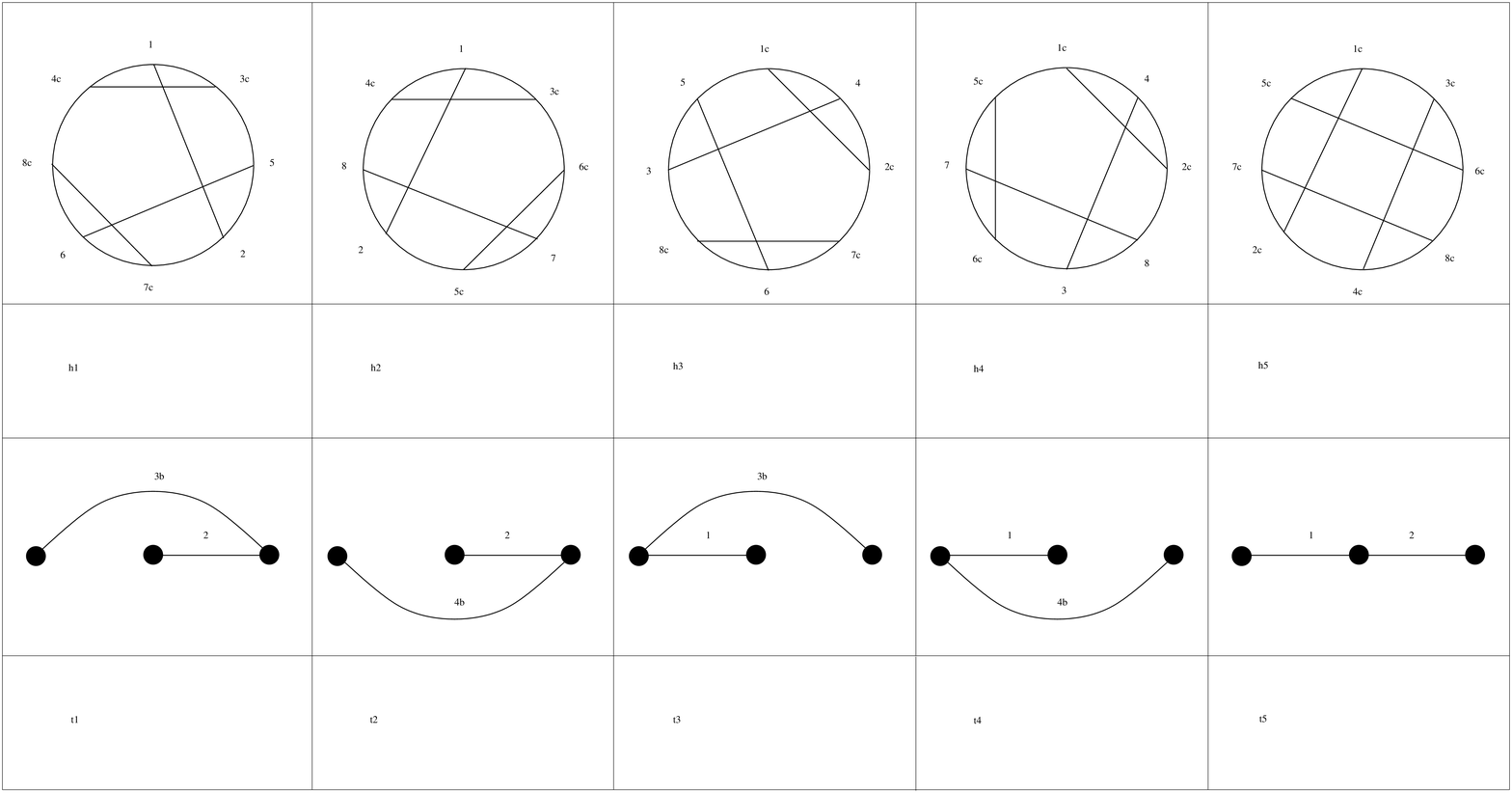}
 \end{center}

\bibliography{dkh}
\bibliographystyle{plain}

\end{document}